\documentclass[11pt,reqno,tbtags]{amsart}
\usepackage{amssymb}
\usepackage{natbib}
\bibpunct[, ]{[}{]}{;}{n}{,}{,}

\numberwithin{equation}{section}

\allowdisplaybreaks

\newtheorem{theorem}{Theorem}[section]
\newtheorem{lemma}[theorem]{Lemma}

\theoremstyle{definition}

\newtheorem{remark}[theorem]{Remark}

\newtheorem*{ack}{Acknowledgement}

\theoremstyle{remark}

\newenvironment{romenumerate}{\begin{enumerate}% gives (i), (ii) etc.
 }{\end{enumerate}}

\newcounter{oldenumi}
{\setcounter{oldenumi}{\value{enumi}}
\begin{romenumerate} \setcounter{enumi}{\value{oldenumi}}}
{\end{romenumerate}}

\newcounter{thmenumerate}
\newenvironment{thmenumerate}
{\setcounter{thmenumerate}{0}%
 \def\item{\par% \ifnum\thethmenumerate=0\else\par\fi %\noindent\fi
 \refstepcounter{thmenumerate}\textup{(\roman{thmenumerate})\enspace}}
}
{}

\newcounter{xenumerate}   %no left indentation; thus wider lines

\newcommand{\refT}[1]{Theorem~\ref{#1}}

\newcommand{\refL}[1]{Lemma~\ref{#1}}
\newcommand{\refR}[1]{Remark~\ref{#1}}
\newcommand{\refS}[1]{Section~\ref{#1}}

\newcommand{\refand}[2]{\ref{#1} and~\ref{#2}}

\begingroup
  \divide\count255 by 60
  \multiply\count255 by -60
  \advance\count255 by \time
  \ifnum \count255 < 10 \xdef\klockan{\the\count1.0\the\count255}
  \else\xdef\klockan{\the\count1.\the\count255}\fi
\endgroup

\newcommand\set[1]{\ensuremath{\{#1\}}}

\newcommand\xpar[1]{(#1)}
\newcommand\bigpar[1]{\bigl(#1\bigr)}
\newcommand\Bigpar[1]{\Bigl(#1\Bigr)}

\newcommand\lrpar[1]{\left(#1\right)}

\newcommand\abs[1]{|#1|}

\def\rompar(#1){\textup(#1\textup)}    % usage: \rompar(...)
\newcommand\xfrac[2]{#1/#2}

\newcommand\parfrac[2]{\lrpar{\frac{#1}{#2}}}

\newcommand\Bigparfrac[2]{\Bigpar{\frac{#1}{#2}}}

\def\xexp(#1){e^{#1}}

\newcommand\floor[1]{\lfloor#1\rfloor}

\newcommand\ntoo{\ensuremath{{n\to\infty}}}

\newcommand\norm[1]{\|#1\|}

\newcommand\downto{\searrow}
\newcommand\upto{\nearrow}

\newcommand\ie{i.e.\spacefactor=1000}
\newcommand\eg{e.g.\spacefactor=1000}

\newcommand{\aex}{a.e.\spacefactor=1000}
\newcommand\whp{{w.h.p.\spacefactor=1000}}

\newcommand{\tend}{\longrightarrow}
\newcommand\dto{\overset{\mathrm{d}}{\tend}}
\newcommand\pto{\overset{\mathrm{p}}{\tend}}

\newcommand\op{o_{\mathrm p}}

\newcommand\bbR{\mathbb R}

\newcommand\bbZ{\mathbb Z}

\newcounter{CC}
\newcounter{cc}
\newcommand\E{\operatorname{\mathbb E{}}}
\renewcommand\P{\operatorname{\mathbb P{}}}

\newcommand\Po{\operatorname{Po}}

\newcommand\ga{\alpha}
\newcommand\gb{\beta}
\newcommand\gd{\delta}

\newcommand\gf{\varphi}
\newcommand\gam{\gamma}

\newcommand\gs{\sigma}

\newcommand\eps{\varepsilon}

\newcommand\cC{\mathcal C}

\newcommand\cE{\mathcal E}
\newcommand\cF{\mathcal F}

\newcommand\gdx{\Delta^*}
\newcommand\ett[1]{\boldsymbol1[#1]} 

\def\[#1]{[\![#1]\!]}

\newcommand\qq{^{1/2}}
\newcommand\qqq{^{1/3}}

\newcommand\qqqd{^{4/3}}

\newcommand\qw{^{-1}}
\newcommand\qww{^{-2}}

\renewcommand{\=}{:=}

\newcommand\intot{\int_0^t}

\newcommand\oi{[0,1]}
\newcommand\ot{[0,t]}
\newcommand\ooo{[0,\infty)}

\newcommand\dd{\,\textup{d}}

\newcommand\lhs{left-hand side}
\newcommand\rhs{right-hand side}

\newcommand\gnp{\ensuremath{G(n,p)}}
\newcommand\gnm{\ensuremath{G(n,m)}}

\newcommand\gx{G^*}
\newcommand\gy{\tilde G}
\newcommand\kk{\kappa}
\newcommand\knk{\komp_k(F_n)}

\newcommand\bfx{\mathbf{x}}

\newcommand\xs{\ensuremath{{\bfx}_n}}

\newcommand\xss{\ensuremath{(\xs)_{n\ge 1}}}

\newcommand\vxs{\ensuremath{\mathcal V}}% vertex space
\newcommand\bbzp{\bbZ_+}
\newcommand\bbzpp{{\bbZ_+}}
\newcommand\hmu{\widehat\mu}
\newcommand\hnu{\widehat\nu}
\newcommand\hnun{\hnu_n}
\newcommand\nun{\nu_n}
\newcommand\nn{^{(n)}}
\newcommand\pij{p_{ij}}
\newcommand\gynt{\gy(n,t;F_n)}
\newcommand\komp{\upsilon}
\newcommand\sssx{s}
\newcommand\sss[1]{\sssx_{#1}}
\newcommand\bsssx{\bar s}
\newcommand\bsss[1]{\bsssx_{#1}}
\newcommand\bs{\bar s}
\newcommand\bx{\bar x_1}

\newcommand\rhokk{\rho_{\kk}}
\newcommand\otc{[0,\tc)}
\newcommand\otf{[0,t_f)}
\newcommand\BFp{Bohman--Frieze process}
\newcommand\BFq{Bohman--Frieze}
\newcommand\BF{\mathsf{BF}}
\newcommand\BFx[1]{\BF\xpar{#1}}
\newcommand\BFt{\BFx t}
\newcommand\tBF{\widetilde{\BF}}
\newcommand\tbfd{\tBF(\tc+\gd)}
\newcommand\fbf{\rho_{\BF}}
\newcommand\BFG{H}
\newcommand\tc{t_{\mathsf c}}
\newcommand\ef{e(F)}
\newcommand\XX{X^*}
\newcommand\citetq[2]{\citeauthor{#2} \cite[{\frenchspacing #1}]{#2}} 
 
\newcommand\ff{\Phi}
\newcommand\gdn{\gd_n}
\newcommand\C{K}
\newcommand\cg{a}
\newcommand{\Holder}{H\"older}

\newcommand\ER{Erd\H os--R\'enyi}

\newcommand\REM[1]{{\raggedright\texttt{[#1]}\par\marginal{XXX}}}

\newenvironment{comment}{\setbox0=\vbox\bgroup}{\egroup} %deletes!

\newcommand\urladdrx[1]{{\urladdr{\def~{{\tiny$\sim$}}#1}}}

\begin{document}
\title
{Phase transitions for modified Erd\H{o}s--R\'enyi processes}

\date{May 25, 2010} % (typeset \today{} \klockan)} %; revised ...

\author{Svante Janson}
\address{Department of Mathematics, Uppsala University, PO Box 480,
SE-751~06 Uppsala, Sweden}
\email{svante.janson@math.uu.se}
\urladdrx{http://www.math.uu.se/~svante/}

\author{Joel Spencer}
\address{Joel Spencer,
Courant Institute,
251 Mercer St.,
New York, NY 10012, USA}
\email{spencer@cs.nyu.edu}
\urladdr{http://www.cs.nyu.edu/cs/faculty/spencer/}
\thanks{This research was mainly done while the authors visited Institute
  Mittag-Leffler, Djursholm, Sweden, 2009} 

%\keywords{<keywords>}

\begin{abstract} 
A fundamental and very well studied region of the Erd\H{o}s-R\'enyi
process is the phase transition at $m\sim \frac{n}{2}$ edges in which a giant 
component suddenly appears.  We examine the process beginning with an initial 
graph.  We further examine the Bohman--Frieze process in which edges between 
isolated vertices are more likely.  While the positions of the phase
transitions  
vary, the three processes belong, roughly speaking, to the same universality
class.   
In particular, the growth of the giant component in the barely supercritical
region  
is linear in all cases. 
\end{abstract}

\maketitle

\section{Introduction}\label{S:intro}

The standard \ER{} process $(\gnm)_{m=0}^{\binom n2}$ starts with an
empty graph $G(n,0)=E_n$ with $n$ vertices
and adds edges one by one in random order,
uniformly over all possibilities, \ie, drawing the edges uniformly
without replacement. (Hence, $\gnm$ has $n$ vertices and $m$ edges.)
This random graph model has been studied a great deal, starting with
\citet{ER59,ER60}, see for example the monographs by \citet{Bollobas}
and \citet{JLR}.

The purpose of this paper is to study two modifications of this process.
We are interested in the sizes (orders) of the components of
the random graphs; in particular whether there exists a giant component of size
comparable to the entire graph and, if so, how large it is. (We ignore the
internal structure of the 
components.)  
We denote
the components of a graph $G$ by $\cC_i(G)$,  $i=1,\dots,{\komp(G)}$, 
where thus $\komp(G)$ is the number of components of $G$, 
and their sizes by $C_i(G)\=|\cC_i(G)|$, $1\le i\le\komp(G)$;
we will always assume that the components are ordered such that $C_1\ge
C_2\ge\dots$.
(For convenience we also define $C_i(G)=0$ when $i>\komp(G)$.) 
We will often, as just done, omit the argument $G$ when the graph is clear from
the context. We further denote the edge set of $G$ by $E(G)$,  the number
of edges by $e(G)\=|E(G)|$, and the number of vertices by $|G|$ (the
\emph{order} or \emph{size} of $G$).

We recall the fundamental result for \gnm{} \cite{ER60} 
that if \ntoo{} and $m\sim cn/2$ for some constant $c$, 
then $C_1=\rho(c)n+\op(n)$, where $\rho(c)=0$ if $c\le1$,
%(in this case thus $C_1=\op(n)$), 
and $\rho(c)>0$ if $c>1$. (Furthermore, $C_2=\op(n)$ for every $c$.)
This is usually expressed by saying that there is a \emph{threshold}
or \emph{phase transition} at $m=n/2$. See further
\cite{ER60,Bollobas,JLR}. 
Moreover, as $\gd\downto0$, $\rho(1+\gd)\sim2\gd$ 
(see \cite[Theorem  3.17]{SJ178} for a generalization to certain other
random graphs). 
(For the notation $\op(n)$, and other standard notations used below such as
\whp, see \eg{} \cite{JLR} and \cite{SJN6}.)

In the first modification of the \ER{} process, we assume that some
(non-random) edges are 
present initially; additional edges then are added randomly as above.
We actually consider three slightly different versions of this process; see
\refS{SER} for details. Our main result for these processes (\refT{T1})
characterizes the existence and size of a giant component in terms of the
initial edges (more precisely, the sizes of the components defined by them)
and the number of added random edges.
We define the \emph{susceptibility} $s_2$ as the average size of the
component containing a random vertex in the initial graph, see
\eqref{Sk}--\eqref{sk2}, and show the existence
of a threshold when $\tc n/2$ edges are added, where  $\tc\=\sss2\qw$.
(This was also done, under a technical assumption, in \citet{SW}.)
% \cite[Theorem 3.1]{SW}.
Moreover, we give upper and lower bounds for the size of the giant
component after the threshold in terms of $s_2$ and two related quantities
(higher moments of the component size) $s_3$ and $s_4$ for the initial
graph, also defined in \eqref{Sk}--\eqref{sk2}.   

Our second modification is known as the Bohman--Frieze process, after
\citet{BF}.  The initial graph 
on $n$ vertices is empty.  At each round two edges $e_1=\{v_1,w_1\}$ and
$e_2=\{v_2,w_2\}$ are selected independently and uniformly.  \emph{If} both
$v_1$ and $w_1$ are isolated vertices the edge $e_1$ is added to the graph;
{\em otherwise} the edge $e_2$ is added to the graph.  We let $\BF_m$ denote
this process when $m$ edges are added.  This is a natural example of an
Achlioptas process, in which a choice may be made from two randomly chosen
potential edges.  In \citet{BF} and \citet{BFW} it was shown that the phase
transition is deferred beyond $m\sim \xfrac{n}{2}$.  
More precisely,
it is proved in \citet{SW} that the Bohman--Frieze process has a phase
transition at some $\tc\approx 1.1763$.
In the present paper we study further what happens just after the phase
transition, using the result just described for the \ER{} process with
initial edges. The idea is, as in \cite{SW},
that to study the process at a time $t_1>\tc$, we
stop the process at a suitable time $t_0$ just \emph{before} the phase
transition, and 
then approximate the evolution between $t_0$ and $t_1$ by an \ER{} process,
using the graph obtained at time $t_0$ as our initial graph.
In order to apply \refT{T1}, we then need information on $s_2$, $s_3$ and
$s_4$ in the subcritical phase.
The analysis in \citet{SW} of the 
Bohman--Frieze process (and a class of generalizations of it) 
is based on studying the susceptibility $s_2$ in the subcritical region.
We will use some results from \cite{SW}, reviewed in \refS{SBF}, and extend
them to $s_3$ and $s_4$ in order to obtain the required results needed to
apply \refT{T1}.

In particular, we show that after the phase transition, the giant component
grows at a linear rate, just as for the 
\ER{} process. The precise statement is given by Theorem \ref{TBF}.
The original Erd\H{o}s--R\'enyi process, the process from an appropriate
starting point, and the Bohman--Frieze process appear to be in what
mathematical physicists loosely call the same universality class.  While
the placement of the phase transitions differ the nature of the phase
transitions appears to be basically the same.  A very different picture
was given for a related process
in \cite{SSA}.  There, as in the Bohman--Frieze process, two random
potential edges $e_1=\{v_1,w_1\}$ and $e_2=\{v_2,w_2\}$ are given.  However
the edge is selected by the \emph{Product Rule}: 
we select that edge for which the
product of the component sizes of the two vertices is largest.  Strong
computational evidence is presented indicating clearly that this process is
not in the same univerality class as the three processes we compare.  We feel,
nonetheless, that there is likely to be a wide variety of processes in the
same universality class as the bedrock Erd\H{o}s--R\'enyi process.

The main results are stated in Sections \ref{SER} and \ref{SBF}, and proved
in Sections \refand{SpfT1}{SpfBF}.

Our results are asymptotic, as the size grows. All unspecified limits are as
$\ntoo$.

\begin{ack}
This research was mainly done at Institute Mittag-Leffler, Djursholm, Sweden,
during the program \emph{Discrete Probability}, 2009.
We thank other participants, in particular Oliver Riordan, for helpful
comments. 

We thank Will Perkins for the numerical calculations in \refR{Rnum}. 
\end{ack}

\section{\ER{} process with an initial graph}\label{SER}

The purpose of this section is to study the \ER{} process when some
edges are present initially. We define three
different but closely related versions of the process.

Let $F$ be a subgraph of $K_n$ with vertex set
$V(F)=V(K_n)=\set{1,\dots,n}$.
Define $(G(m,n;F))_{m=0}^{\binom n2-\ef}$ by starting with $G(n,0;F)\=F$ 
and adding the $\binom n2-\ef$ edges in $E(K_n)\setminus E(F)$ one by one in
random order, \ie, by drawing without replacement.

For our purposes it will be convenient to consider two modifications
of this random graph process. (Both modifications are well-known for \gnm.)
We define $(\gx(n,m;F))_{m=0}^\infty$ by starting with $\gx(n,0;F)\=F$
and then adding at each time step an edge randomly drawn 
(with replacement) from
$E(K_n)$, provided this edge is not already present (in which case
nothing happens). In particular, $\gx(n,m)\=\gx(n,m;E_n)$ is defined
as \gnm{} but drawing the edges with replacement. 
In general, we have $E(\gx(n,m;F))=E(\gx(n,m))\cup E(F)$.

%\begin{remark}
Note that
 the number of
edges in $\gx(n,m)$ may be less than $m$.
Alternatively, we may regard $\gx(n,m;F)$ as a multigraph and add the
edges whether they already are present or not; then the number of
edges is always exactly $m+\ef$. Since we will study the component
sizes only, this makes no difference for the present paper.
%\end{remark}

The second modification is to use continuous time. We may think of the
$\binom n2$ edges as arriving according to independent Poisson
processes with rates $1/n$; thus edges appear at a total rate $\binom
n2/n=\frac{n-1}2$ and each edge is chosen uniformly at random and
independently of all previous choices.
We define $\gy(n,t;F)$ to be $F$ together with all edges that have
arrived in $[0,t]$.
(As above, we can consider either a multigraph version or the
corresponding process of simple graphs, obtained by ignoring all edges
that already appear in the graph.)
Hence, if $i$ and $j$ are two vertices that are not already joined by
an edge in $F$, then the probability that they are joined in
$\gy(n,t;F)$ is $1-e^{-t/n}=t/n+O(t^2/n^2)$, and these events are
independent for different pairs $i,j$. (Starting with the empty graph
we thus obtain $\gnp$ with $p=1-e^{-t/n}$. We could change the time
scale slightly to obtain exactly $G(n,t/n)$, and asymptotically we
obtain the same results for the two versions.)

Note that if $N(t)$ is the total number of edges arriving in $\ot$,
then $N(t)\sim \Po\bigpar{\binom n2 t/n}=\Po\bigpar{\frac{n-1}2t}$,
and, with an obvious coupling of the processes, $\gy(n,t;F)=\gx(n,N(t);F)$.
For constant $t$, $N(t)/(n/2)\pto t$ as \ntoo{} by the law of large
numbers. 
Moreover, the expected
number of repeated edges in $\gx(n,m;F)$ is at most $\binom m2/\binom
n2+m\abs{E(F)}/\binom n2$; if for example, as in \refT{T1}
%our application in \refS{SBF}, 
$m=O(n)$ and $\abs{E(F)}=O(n)$, then this is $O(1)$, which
will be negligible. Standard arguments, comparing the processes at
times $t$ and $(1\pm\eps)t$, show that for the properties considered
here, and asymptotically as \ntoo, we then obtain the same results for
$G(n,\floor{nt/2};F)$, 
$\gx(n,\floor{nt/2};F)$, 
and $\gy(n,t;F)$.

We define, for a graph $G$ with components of sizes $C_1,\dots,C_\komp$, 
and $k\ge1$,
\begin{equation}\label{Sk}
S_k=  S_k(G)\=\sum_i C_i^k,
\end{equation}
summing over all components of $G$.
Thus $S_1(G)=|G|$, the number of vertices. 
We normalize these sums by dividing by $|G|$ and define
\begin{equation}\label{sk}
  \sssx_k=\sssx_k(G)\=\frac{S_k(G)}{|G|}=\frac{S_k(G)}{S_1(G)}.
\end{equation}
Hence, $\sss1(G)=1$ for every $G$. 
Note that 
\begin{equation}\label{sk2}
\sss k(G)=\sum_i \frac{C_i}{|G|}C_i^{k-1},   
\end{equation}
which is
the $(k-1)$:th
moment of the size of the component containing a randomly chosen vertex.
In particular, $\sss2(G)$ is the average size of the component containing a
random vertex.
The number $\sss2(G)$ is called
the \emph{susceptibility}; see \eg{} \cite{SJ218,SJ232,SJ241} for results on
the susceptibility in $\gnm$ and some other random graphs. 

It follows from the definitions \eqref{Sk} and
\eqref{sk} that $S_k$ and $s_k$ are (weakly) increasing in
$k$; in particular, $\sss{k}(G)\ge\sss1(G)=1$ for
every $k$ and $G$. Moreover, \Holder's inequality and \eqref{sk2}
imply that the stronger result that $\sssx_k^{1/k}$ 
(and even $\sssx_k^{1/(k-1)}$, $k\ge2$)
is
(weakly) increasing in $k$.

Note further that the number of edges in a component of size
$C_i$ is at most $\binom{C_i}2\le C_i^2$;
hence, for any graph $G$, 
\begin{equation}\label{eg}
  |E(G)|\le S_2(G).
\end{equation}

We will use these functionals for the initial graph $F$ to characterize the
existence and size of a giant component in the random graph processes
starting with $F$. An informal summary of the following theorem (our main
result in this section) is that there is a phase transition at
$\tc\=1/s_2(F)$,   
and that for $t=\tc+\gd$ with $\gd$ small, there is a giant component of
size $\approx 2(s_2(F)^3/s_3(F))\gd n$. For the special case when $F=E_n$ is
empty, $s_2=s_3=1$ and we recover the well-known result for the \ER{}
process mentioned above that there is a phase transition at $\tc=1$ (\ie, at
$n/2$ 
edges) and further for $t=1+\gd$, there is a giant component of size
$\approx 2\gd n$. 
The formal statement is asymptotic, and we thus consider a
sequence $F_n$.

\begin{theorem}\label{T1}
  Suppose that for each $n$ (at least in some subsequence), $F_n$
  is a given graph with $n$ vertices, and 
suppose that $\sup_n\sss3(F_n)<\infty$. Let the random variable
$Z_n$ be the size of the   component containing a random vertex in $F_n$. 

Consider the random graph
  processes $\gy(n,t;F_n)$.
Then, for any fixed $t>0$, the following hold as \ntoo, 
with $s_k\=s_k(F_n)$,
\begin{romenumerate}
  \item\label{t1sub}
If $t\le1/\sss2$, then $C_1(\gy(n,t;F_n))=\op(n)$.
\item\label{t1gan}
If\/ $t>1/s_2$, then there is
a unique $\rho_n>0$
such that 
\begin{equation*}
  \rho_n=1-\E e^{-\rho_n t Z_n},
\end{equation*}
and we have
\begin{equation*}
  C_1(\gy(n,t;F_n))=\rho_nn+\op(n).
\end{equation*}
  \item\label{t1+-p}
If $t> 1/\sss2$, let $\gd_n\=t-1/\sss2>0$.
Then 
\begin{equation*}
\frac{ C_1(\gy(n,t;F_n))}n\ge
2\gd_n 
\frac{\sss2^3}{\sss3}
\lrpar{1-2\gd_n\sss2}
+\op(1).
\end{equation*}
If further 
%$\gd_n\le\frac38 \sss3^2/(\sss2^2\sss4)$,
$\gd_n\sss2^2\sss4/\sss3^2\le\frac38$,
then also
\begin{equation*}
\frac{ C_1(\gy(n,t;F_n))}n
\le 2\gd_n 
\frac{\sss2^3}{\sss3} \lrpar{1+\frac83\gd_n\frac{\sss2^2\sss4}{\sss3^2}}
+\op(1).
\end{equation*}
\item\label{t1+-whp}
In \ref{t1+-p}, if in addition $\liminf_\ntoo\gd_n>0$, then 
moreover
\whp
\begin{equation*}
\frac{ C_1(\gy(n,t;F_n))}n\ge
2\gd_n 
\frac{\sss2^3}{\sss3}
\lrpar{1-2\gd_n\sss2}
\end{equation*}
and, if  
$\gd_n\sss2^2\sss4/\sss3^2\le\frac38$,
\begin{equation*}
%2\gd_n \frac{\sss2^3}{\sss3}\lrpar{1-2\gd_n\sss2}\le
\frac{ C_1(\gy(n,t;F_n))}n
\le 2\gd_n 
\frac{\sss2^3}{\sss3} \lrpar{1+\frac83\gd_n\frac{\sss2^2\sss4}{\sss3^2}}
.
\end{equation*}
\end{romenumerate}

The same results hold for the random graph processes
$G(n,\floor{nt/2};F)$ and 
$\gx(n,\floor{nt/2};F)$.
\end{theorem}
The proof is given in \refS{SpfT1}.
Note that  by \eqref{sk2},
\begin{equation}
  \label{ezk}
\E Z_n^k=s_{k+1}(F_n), \qquad k\ge1.
\end{equation}

\section{The Bohman--Frieze process}\label{SBF}

Recall the definition of the \BFp{} from \refS{S:intro}, see \cite{BF,BFW,SW}:
we are at each round presented with two random edges 
$e_1=\{v_1,w_1\}$ and $e_2=\{v_2,w_2\}$ in the complete graph $K_n$
and choose one of them; we choose
$e_1$ if both its endpoints $v_1$ and $w_1$ are isolated, and otherwise we
choose $e_2$.
We let $\BF_m$ denote the random graph created by this process
when $m$ edges are added.
%; for convenience we let $\BF_s\=\BF_{\floor s}$ for any real $s\ge0$.
(The size $n$ is not shown explicitly.)
We further define, using the natural time scale,
$\BFt\=\BF_{\floor{nt/2}}$. (For convenience, we sometimes omit
rounding to integers in expressions below.)

Note that if we add $e_1$, then it always joins two previously isolated
vertices, while if we add $e_2$, it is uniformly distributed and
independent of the existing graph. We call the added edges $e_2$
\emph{\ER{} edges}, since all edges in the \ER{} process are of this type.

\begin{remark}
  \label{Rloops}
We have talked about edges 
$e_1$ and $e_2$, but it is technically convenient in the proofs
to allow also loops (as in \cite{SW}); we thus assume in the proofs below
that in each round,
the vertices $v_1,w_1,v_2,w_2$ are independent, uniformly distributed,
random vertices. 
It is easily seen that the results proved for this version hold also if we
assume that there are no loops, for example by conditioning on the event
that no loops are presented during the first $nt/2$ rounds; we omit the
details. 
\end{remark}

For a graph $G$, let $n_i=n_i(G)$ be the number of vertices in components of
order $i$, and let $x_i=x_i(G)\=n_i(G)/|G|$ be the proportion of the total
number of vertices that are in such components. 
(Thus, $s_k(G)=\sum_i i^{k-1}x_i(G)$.) 
For the Bohman--Frieze process, we  need only $n_1$, the number of
isolated vertices, and the corresponding proportion $x_1\=n_1/n$.

For the Bohman--Frieze process (and some generalizations of it), it is shown
in \citet{SW} that the random variables
$x_1(\BFt)$ (for any fixed $t<\infty$) and $s_2(\BFt)$ 
(for any fixed $t<\tc$) converge in probability, as \ntoo,
to some deterministic values $\bx(t)$ and $\bs_2(t)$; these limit values are
given as solutions to differential equations.
We extend this to $s_3$ and $s_4$ as follows. 

We first define, as in \cite{SW}, 
the deterministic function $\bx(t)$ as the solution to the
differential equation
\begin{equation}
  \label{x1}
\bx'(t)=-\bx^2(t)-\bigpar{1-\bx^2(t)}\bx(t),
\qquad t\ge0,
\end{equation}
with initial condition $\bx(0)=1$; by \cite[Theorem 2.1]{SW}, $\bx(t)$ is
defined and  positive for all $t\ge0$, and by \cite[Theorem  1.1]{SW}, 
$x_1(\BFt)\pto \bx(t)$ for every fixed $t\ge0$.

We further define functions $\bs_2(t)$, $\bs_3(t)$, $\bs_4(t)$ as the
solutions to the differential equations
\begin{align}
  \bs_2'(t)&=\bx^2(t)+\bigpar{1-\bx^2(t)}\bs_2^2(t), \label{s2}\\
  \bs_3'(t)&=3\bx^2(t)+3\bigpar{1-\bx^2(t)}\bs_2(t)\bs_3(t), \label{s3}\\ 
 \bs_4'(t)&=7\bx^2(t)+\bigpar{1-\bx^2(t)}\bigpar{4\bs_2(t)\bs_4(t)+3\bs_3^2(t)}, 
\label{s4}
\end{align}
with initial conditions 
\begin{equation}\label{s2340}
 \bs_2(0)=\bsss3(0)=\bsss4(0)=1. 
\end{equation}
The function $\bs_2(t)$ 
is studied in \citet[Theorem 2.2]{SW}, 
and it is shown there that 
it explodes at some finite $\tc$, \ie, the solution $\bs_2(t)$ is (uniquely)
defined for $t\in[0,\tc)$, but $\bs_2(t)\upto+\infty$ as $t\upto\tc$; it is
further shown \cite[Theorem 1.1]{SW} that this $\tc$ is the time of the
phase transition for the \BFp, when a giant component first appears, and
that for any fixed $t<\tc$, $\bs_2(\BFt)\pto \bs_2(t)$.
We extend these results to $\bs_3$ and $\bs_4$ as follows.

\begin{theorem}\label{TSa}
  The functions $\bs_2(t)$, $\bs_3(t)$, $\bs_4(t)$ are uniquely defined by
  \eqref{s2}--\eqref{s2340} for all $t\in[0,\tc)$.
As $t\upto\tc$, there exist positive constants $\ga$ and $\gb$ such that
\begin{align*}
  \bs_2(t)&\sim \frac{\ga}{\tc-t}, \\
  \bs_3(t)&\sim \gb \bs_2(t)^3\sim \frac{\gb\ga^3}{(\tc-t)^3}, \\
  \bs_4(t)&\sim 3\gb^2 \bs_2(t)^5\sim \frac{3\gb^2\ga^5}{(\tc-t)^5}.
\end{align*}
More precisely, $\bs_k(t)=\cg_k(\tc-t)^{-(2k-3)}(1+O(\tc-t))$ for $k=2,3,4$
with 
$\cg_2=\ga$, $\cg_3=\gb\ga^3$, $\cg_4=3\gb^2\ga^5$.

We have $\ga=\bigpar{1-\bx^2(\tc)}\qw$, while $\gb=g(\tc)$ is given by 
\eqref{g'} and \eqref{g}. %,  \eqref{G} and \eqref{f'}. 
\end{theorem}

\begin{theorem}\label{TSb}
For any fixed $t\in\otc$, and $k=2,3,4$, $s_k(\BFt)\pto \bs_k(t)$.  
\end{theorem}

\begin{remark}
  \label{RSb}
It is straightforward to extend \refT{TSb} to any $k\ge2$, with $\bs_k(t)$
given by a differential equation similar to \eqref{s2}--\eqref{s4}
(involving $\bs_j$ for $j<k$, so the functions are defined recursively).
We leave the details to the reader since we only use $k\le4$ in the present
paper. 
\end{remark}

Proofs are given in \refS{SpfBF}.
Using these results for the subcritical phase, we obtain the following for
the supercritical phase; 
again the proof is given in \refS{SpfBF}.

\begin{theorem}
  \label{TBF}
There exists  constants $\gam=2 (1-\bx^2(\tc))/\gb>0$ and $\C<\infty$ such
that for 
any fixed $\gd>0$, \whp{} 
\begin{equation*}
\gam \gd -\C\gd\qqqd
\le 
 \frac{ C_1(\BFx{\tc+\gd})}{n} \le \gam \gd +\C\gd\qqqd.
\end{equation*}
\end{theorem}

\begin{remark}\label{Rnum}
  Numerical calculations of Will Perkins give 
$\tc\approx 1.1763$, 
$\bx(\tc)\approx 0.2438$,  
$\ga \approx 1.063$, 
$\gb \approx 0.764$, %0.7637?
$\cg_2=\ga$,
$\cg_3\approx 0.917$, 
$\cg_4\approx 2.375$
%$\bs_3(t)\approx 0.917(\tc-t)^{-3}$ 
and $\gam \approx 2.463$.  %2.46318.
\end{remark}

There is an obvious conjecture (made explicit in \cite{SW}) that 
$\xfrac{ C_1(\BFx{t})}{n} \pto \fbf(t)$ for some function $\fbf:\ooo\to\oi$;
equivalently, $C_1(\BFx{t})=\fbf(t)n+\op(n)$.
(For $t<\tc$, clearly this holds with $\fbf(t)=0$.)
In \citet{SW} it was further
conjectured that $\lim_{\gd\rightarrow \tc^+} \fbf(t) = 0$;
in the language of Mathematical Physics, this says that the phase 
transition is not first order.  
If such an $\fbf$ exists, \refT{TBF} resolves the latter conjecture
positively and further gives the asymptotic behavior 
$\fbf(\tc+\gd)\sim\gam\gd$ as $\gd\rightarrow 0^+$. 
 
\begin{remark}
  \label{RBF}
We further conjecture that the function $\fbf$ is smooth
on $[\tc,\infty)$; if this is the case, then \refT{TBF} shows
that  $\fbf'(\tc^+)=\gam$. 
This conjecture
would imply that $\gd\qqqd$ in  \refT{TBF} could be replaced by $\gd^2$;
unfortunately, our approximations are not sharp enough to show this.
\end{remark}

\section{Proof of \refT{T1}}\label{SpfT1}

We begin with a simple lemma (related to results in
\cite[Section 5]{SJ178}).

\begin{lemma}\label{L1}
  Let\/ $Y\ge0$ be a random variable with $1<\E Y\le\infty$.
  \begin{thmenumerate}
\item\label{l1ga}	
There is a unique $\rho>0$ such that 
\begin{equation}\label{gal}
  \rho=1-\E e^{-\rho Y}.
\end{equation}
\item\label{l1-}
If\/ $\E Y^2<\infty$, then
\begin{equation*}
  \rho>\frac{2 (\E Y-1)}{\E Y^2} .
\end{equation*}
\item\label{l1+}
If\/ $\E Y^3<\infty$ and 
$8(\E Y-1)\E Y^3\le 3(\E Y^2)^2$,
then
\begin{equation*}
  \begin{split}
  \rho&<\frac{3\E Y^2-\sqrt{9(\E Y^2)^2-24(\E Y-1)\E Y^3}}{2\E Y^3}	
%\\&
%=\frac{12(\E Y-1)}{3\E Y^2+\sqrt{9(\E Y^2)^2-24(\E Y-1)\E Y^3}}
\\&=
\frac{4(\E Y-1)}{\E Y^2+\sqrt{(\E Y^2)^2-\frac83(\E Y-1)\E Y^3}}
\\&\le
\frac{2(\E Y-1)}{\E Y^2}\lrpar{1+\frac{8(\E Y-1)\E Y^3}{3(\E Y^2)^2}}
.
  \end{split}
\end{equation*}

\item\label{l1lim}
Let\/ $Y_n$, $n\ge1$, be random variables with $Y_n\ge0$ and $\E Y_n>1$ and 
let $\rho_n>0$ be the corresponding numbers such that 
$\rho_n=1-\E e^{-\rho_n  Y_n}$. 
If\/ $Y_n\dto Y$ for some $Y$ with $\E Y>1$, then $\rho_n\to\rho>0$ satisfying
\eqref{gal}. 
On the other hand, if\/ $Y_n\dto Y$ with $\E Y\le 1$, then $\rho_n\to0$.
  \end{thmenumerate}
\end{lemma}

\begin{remark}
  In fact, \eqref{gal} is the standard equation for the survival probability
  of a Galton--Watson process with a mixed Poisson $\Po(Y)$ offspring
  distribution. Parts (i) and (iv)
  follow easily from standard results on branching processes. We prefer,
  however, to give direct proofs (also easy).  
Note further that $\rho=0$ always is another solution to \eqref{gal}. If $\E
Y\le1$, 
then $\rho=0$ is the only non-negative solution, either by branching process
theory, or because
\begin{equation}
  \label{ga0}
1- \E e^{-s Y}
= \E(1- e^{-s Y}) \le \E(sY) \le s
\end{equation}
for every $s\ge0$, with strict inequality unless $sY=0$ \aex{} and
$\E(sY)=s$, which together imply $s=0$.
\end{remark}

\begin{proof}
The function
$\gf(s)\=1-\E e^{-sY}$, $s\in[0,\infty)$, is increasing and concave with
  $0\le\gf(s)<1$, $\gf(0)=0$ and $\gf'(0)=\E Y>1$. Consequently, $\gf(s)>s$
  for small $s>0$, but $\gf(s)<s$ for $s>1$, say, and there is a unique
  $\rho>0$ such that $\gf(\rho)=\rho$. This proves \ref{l1ga}.
Note that $\gf(s)>s$ for $0<s<\rho$ and $\gf(s)<s$ for $s>\rho$.

We next prove \ref{l1lim}. 
If $\E Y>1$, let $0<\eps<\rho$. Then $\gf(\rho-\eps)>\rho-\eps$ and thus,
because $Y_n\dto Y$,
\begin{equation*}
  1-\E e^{-(\rho-\eps)Y_n}
\to   1-\E e^{-(\rho-\eps)Y} >\rho-\eps,
\end{equation*}
so for large $n$, $1-\E e^{-(\rho-\eps)Y_n}>\rho-\eps$ and thus
$\rho-\eps<\rho_n$. 
Similarly, for large $n$, 
$1-\E e^{-(\rho+\eps)Y_n}<\rho+\eps$ and thus
$\rho+\eps>\rho_n$. Since $\eps$ is arbitrarily small, it follows that
$\rho_n\to\rho$.

If instead $\E Y\le1$, then $\gf(s)<s$ for every $s>0$ by \eqref{ga0} and
the comment after it. Hence the same argument shows that for every $\eps>0$, 
$\rho_n<\eps$ for large $n$; thus $\rho_n\to0$.

To see \ref{l1-}, observe that
$e^{-x}\le 1-x+x^ 2/2$ for $x\ge0$, with strict inequality unless
$x=0$, and thus, when $\E Y^2<\infty$,
\begin{equation*}
  \rho=\E\bigpar{1-e^{-\rho Y}}
> \E\lrpar{\rho Y-\frac{\rho ^2Y^2}2}
= \rho \E Y -\frac{\rho ^2}2\E Y^2.
\end{equation*}
Hence, $1>\E Y- \rho\E Y^2/2$, which yields \ref{l1-}.

For \ref{l1+}, we first note that, similarly,
$e^{-x}\ge 1-x+x^ 2/2-x^3/6$ for $x\ge0$, again with strict inequality unless
$x=0$, and thus, provided $\E Y^3<\infty$,
\begin{equation*}
  \rho=\E\bigpar{1-e^{-\rho Y}}
< \E\lrpar{\rho Y-\frac{\rho ^2Y^2}2+\frac{\rho ^3Y^3}6}
= \rho \E Y -\frac{\rho^2}2\E Y^2+\frac{\rho^3}6\E Y^3.
\end{equation*}
This can be written
%\begin{equation*}
%\frac{\E Y^3}6\rho^2-\frac{\E Y^2}2\rho+\E Y-1>0.
%\end{equation*}
\begin{equation}\label{sofie}
\E Y^3\,\rho^2-3{\E Y^2}\,\rho+6(\E Y-1)>0.
\end{equation}
As long as the discriminant 
$9(\E Y^2)^2-24(\E Y-1)\E Y^3\ge0$,
the corresponding quadratic equation (with equality instead of $>$) has
two roots
\begin{equation*}
  \rho_\pm=\frac{3\E Y^2\pm\sqrt{9(\E Y^2)^2-24(\E Y-1)\E Y^3}}{2\E Y^3}
\end{equation*}
and we have either $\rho<\rho_-$ or $\rho>\rho_+$.
In order to rule out the latter possibility, we consider the random variable
$Y_t\=tY$ for $t_0<t\le 1$, where $t_0=1/\E Y$. Note that for $t_0<t\le 1$,
$\E Y_t>1$ and thus there is an $\rho(t)>0$ such that $\rho(t)=1-\E
e^{-\rho(t)Y_t}$; by \ref{l1lim}, $\rho(t)$ is a continuous function of $t$. 
Further, for $t_0<t\le1$,
\begin{equation*}
  \begin{split}
  9(\E Y_t^2)^2 
&= 9t^4  (\E Y^2)^2 
\ge 24t^4(\E Y-1)\E Y^3
= 24(\E Y_t-t)\E Y_t^3
\\&
\ge 24(\E Y_t-1)\E Y_t^3;	
  \end{split}
\end{equation*}
hence the discriminant is non-negative for each $Y_t$, and there are
corresponding roots $\rho_\pm(t)$. These are continuous functions of $t$
and for each $t\in(t_0,1)$, $\rho(t)<\rho_-(t)$ or $\rho(t)>\rho_+(t)$.
As $t\downto t_0$, $\E Y_t\to1$ and $\rho_+(t)\to3\E Y_{t_0}^2/(2\E
Y_{t_0}^3)>0$ while, by \ref{l1lim} again, $\rho(t)\to0$. Hence,
$\rho(t)<\rho_+(t)$ for $t$ close to $t_0$, and by continuity,
$\rho(t)<\rho_+(t)$ for all $t\in(t_0,1]$ (since equality is impossible by
\eqref{sofie}). Consequently, $\rho<\rho_+$ and thus $\rho<\rho_-$. 

Finally, we use straightforward algebra and
the fact that for $x\in\oi$, $\sqrt{1-x}\ge(1-x)/(1+x)$ and
thus
\begin{equation*}
  \frac{1}{1+\sqrt{1-x}}\le\frac{1+x}2.
\qedhere
\end{equation*}
\end{proof}

\begin{proof}[Proof of \refT{T1}]
Note that the assumptions and \eqref{eg} imply that
\begin{equation*}
  |E(F_n)|\le S_3(F_n)=n\sss3(F_n)=O(n). 
\end{equation*}
Hence,
  by the discussion in \refS{SER}, it suffices to consider
  $\gy(n,t;F_n)$.

The main idea is that we may collapse each component $\cC_i(F_n)$ of $F_n$
to a ``supervertex'' with \emph{weight} 
\begin{equation}\label{xi}
x_i=x_i\nn\=|\cC_i(F_n)|=C_i(F_n).  
\end{equation}
The probability of an edge between $\cC_i(F_n)$ and $\cC_j(F_n)$ in
$\gy(n,t;F_n)$ is, for $i\neq j$, 
\begin{equation}
  \label{pij}
\pij(t)=1-e^{-tx_ix_j/n}.
\end{equation}
Hence, to obtain the distribution of component sizes in $\gynt$ we may 
instead consider the random graph $H_n$ with $\komp=\komp(F_n)$ vertices
having weights $x_i$ given by \eqref{xi} and edges added independently with
probabilities $\pij$ given by \eqref{pij}; note that the size of a component
in $\gynt$ is given by the {weight} of the corresponding component in
$H_n$, \ie, the sum of the weights of the vertices in it. 

The random graph $H_n$ is an instance of the general random graph model
studied in \citet{SJ178}; we will use results from \cite{SJ178}, and
therefore we show the relation in some detail.

We will actually consider a subsequence only, for technical reasons, and
thus we at first obtain the result for this subsequence only.
However, this means that if we start with any subsequence of the original
sequence, there exists a subsubsequence where the result holds; this fact
implies that the result actually holds for the full sequence by the
subsubsequence principle, see \eg{} \cite[p.\ 12]{JLR}.
%(This is easily proved by contradiction. If, say, the conclusion of
%(ii) does not hold \whp, then there are $eps>0$,  $\eta>0$ and a subsequence
%where $|C_1-\rho_n n$|>\eps n$ with probability at least $\eta$; (ii) then
%cannot hold for any subsequence of this subsequence.)

We have defined $Z_n$ as the size of the component containing
a random vertex in $F_n$. Let $\hnu_n$ be the
distribution of $Z_n$; thus $\hnu_n$ is the probability
measure on $\bbzp\=\set{1,2,\dots}$ given by 
$\sum_i \frac{C_i}n\gd_{C_i}$. 
By  \eqref{ezk},
 $\E Z_n=\sss2(F_n)\le\sss3(F_n)=O(1)$, which implies that the
%sequence of distributions 
%$\hnu_n$ is tight, see \eg{} \cite[Section 5.8.3]{Gut}. Consequently
%(see \cite[Theorem 5.8.5]{Gut}), we may select a subsequence such that
%$\hnu_n$ converges (weakly) to some probability measure $\hmu$ on $\bbzp$. 
%Equivalently,  $Z_n$ converges in distribution to a random variable $Z$ with
%the distribution $\hmu$.
sequence of random variables
$Z_n$ is tight, see \eg{} \cite[Section 5.8.3]{Gut}. Consequently
(see \cite[Theorem 5.8.5]{Gut}), we may select a subsequence such that
$Z_n$ converges in distribution to some random variable $Z$.
Equivalently, $\hnun$ converges (weakly) to some probability measure $\hmu$
on $\bbzp$, where $\hmu$ is the distribution of $Z$.
Moreover, $\E Z_n^2=s_3(F_n)=O(1)$, and thus \cite[Theorem 5.4.2]{Gut}
$Z_n$  are uniformly integrable; 
consequently \cite[Theorem 5.5.8]{Gut}, 
$s_2(F_n)=\E Z_n\to\E Z$. We denote this limit by $\bs_2$, and have
thus
\begin{align}\label{ez}
s_2(F_n)\to\bs_2=  \E Z .
\end{align}
%while Fatou's lemma for convergence in distribution 
%\cite[Theorem 5.5.8]{Gut} shows that
%\begin{equation}
%  \E Z^2\le\liminf_{\ntoo}\E Z_n^2=\liminf_{\ntoo} s_3(F_n)=\bs_3.
%\end{equation}
%\begin{equation}
%  \E Z^3\le\liminf_{\ntoo}\E Z_n^3=\liminf_{\ntoo} s_4(F_n)=\bs_4.
%\end{equation}

Let $\knk$ be the number of components of order $k$ in $F_n$ and let
$\nu_n$ be the measure on $\bbzp$ defined by
\begin{equation*}
  \nu_n\set k\=\frac{\knk}n.
\end{equation*}
Equivalently, 
$\nu_n\= \frac1n\sum_{i=1}^{\komp}\gd_{C_i}$. 
The total mass of $\nu_n$ is thus $\nu_n(\bbzp)=\komp(F_n)/n\le1$.
(In general, $\nu_n$ is not a probability measure.)

The total size of the components of order $k$ in $F_n$ is $k\knk$, and
thus
\begin{equation*}
  \hnu_n\set k=\P(Z_n=k)=\frac{k\knk}{n}=k\nu_n\set k.
\end{equation*}
Let $\mu$ be the measure on $\bbzp$ given by 
\begin{equation*}
\mu\set k\=\hmu\set k/k,
\qquad k\ge1.  
\end{equation*}
Since we have $\hnun\set k\to\hmu\set k$, we also have 
$$\nun\set k=\hnun\set k/k\to\hmu\set k/k=\mu\set k$$ 
for every $k\ge1$. 
Moreover, if $f:\bbzp\to\bbR$ is any bounded function, and $g(k)\=f(k)/k$,
then the convergence $\hnun\to\hmu$ implies
$$
\int_\bbzpp f(x)\dd\nun(x)
=\int_\bbzpp g(x)\dd\hnun(x)
\to
\int_\bbzpp g(x)\dd\hmu(x)
=\int_\bbzpp f(x)\dd\mu(x).
$$
Hence $\nun\to\mu$ weakly; in particular
\begin{equation}\label{nun}
 \nun(A)\to\mu(A) \qquad \text{for every $A\subseteq\bbzp$}. 
\end{equation}

We let $\xss$ be the sequence $(C_1(F_n),\dots,C_{\komp_n}(F_n))$ of
component sizes of $F_n$, where $\komp_n\=\komp(F_n)$.
We have just shown that the triple $\vxs\=(\bbzp,\mu,\xss)$ is a
\emph{generalized vertex space} in the sense of \cite[p.\ 10]{SJ178};
in particular, the crucial condition \cite[(2.4)]{SJ178} is our \eqref{nun}.

We define the \emph{kernel} $\kk$ on $\bbzp$ by 
\begin{equation}
  \label{kk}
\kk(x,y)\=txy
\end{equation}
(recall that $t$ is fixed); the probability \eqref{pij} of an edge in $H_n$
between (super)vertices with weights $x_i$ and $x_j$ is thus 
$1-\exp(-\kk(x_i,x_j)/n)$, which agrees with \cite[(2.6)]{SJ178}.
Hence, our random graph $H_n$ is the graph denoted $G^\vxs(n,\kk)$ in
\cite{SJ178}. 

We further have, with $x_i=C_i(F_n)$, by \eqref{pij},
\begin{equation*}
  \begin{split}
\frac1n\E e(H_n)
&=
\frac1n\sum_{1\le i<j\le\komp_n}\pij	
=
\frac1n\sum_{1\le i<j\le\komp_n}\bigpar{1-\exp(-tx_ix_j/n)}
\\&
\le
\frac1{n^2}\sum_{1\le i<j\le\komp_n}tx_ix_j
\le
\frac t2\lrpar{\frac1n\sum_{i=1}^{\komp_n}x_i}^2 =\frac t2,
  \end{split}
\end{equation*}
and
\begin{equation}\label{xdmu}
  \int_\bbzpp x\dd\mu(x)
=
  \sum_{x=1}^\infty x\dd\mu\set x
=
  \sum_{x=1}^\infty \dd\hmu\set x
=\hmu(\bbzp)=1
\end{equation}
(since $\hmu$ is a probability measure on $\bbzp$); hence 
\begin{equation}\label{l1}
  \iint_{\bbzp^2}\kk(x,y)\dd\mu(x)\dd\mu(y)=
t\lrpar{\int_\bbzpp x\dd\mu(x)}^2=t
\end{equation}
and
\begin{equation*}
\frac1n\E e(H_n)\le \frac12 \iint_{\bbzp^2}\kk(x,y)\dd\mu(x)\dd\mu(y)
.\end{equation*}
Together with \cite[Lemma 8.1]{SJ178}, this shows that 
$$
\frac1n\E e(H_n)\to \frac12 \iint_{\bbzp^2}\kk(x,y)\dd\mu(x)\dd\mu(y),
$$
and thus, using also \eqref{l1}, the kernel $\kk$ is \emph{graphical}
\cite[Definition 2.7]{SJ178}.

We can now apply the results of \cite{SJ178}.
The kernel $\kk(x,y)$ is of the special type $\psi(x)\psi(y)$ (with
$\psi(x)\=t\qq x$), which is the \emph{rank 1} case studied in \cite[Section
  16.4]{SJ178}, and it follows by \cite[Theorem 3.1 and (16.8)]{SJ178} that
$H_n$ has a giant component if and only if $\norm{T_\kk}>1$, where
$T_\kk$ is the integral operator with kernel $\kk$;
in the rank 1 case $T_\kk$ has the norm,
using
also \eqref{ez},
\begin{equation*}
  \norm{T_\kk}=\int_{\bbzp}\psi(x)^2\dd\mu(x)
=\int_{\bbzp}tx^2\dd\mu(x)
=\int_{\bbzp}tx\dd\hmu(x)
=t\E Z=t\bs_2.
\end{equation*}

Hence there is a phase transition at $\tc\=1/\bs_2$. We consider the cases
$t\le\tc$ and $t>\tc$ separately.

\subsection{The (sub)critical case}
Consider first the case $t\le\bs_2\qw$; then $H_n$ thus has no giant component;
more precisely, 
\begin{equation}\label{c1small}
C_1(H_n)=\op(n).  
\end{equation}
Recall, however, 
that we really are interested in the size of the largest component of
$\gynt$, which is the same as the largest \emph{weight} of a component in
$H_n$. (Note also that the component with largest weight not necessarily is the
component with largest number of vertices.)
Nevertheless, the corresponding estimate follows easily: 
Let $A>0$. Then the
total weight of all vertices in $H_n$ of weight larger than $A$ is
\begin{equation*}
  \begin{split}
	\sum_i x_i\ett{x_i> A}
&=	\sum_{k>A} k\knk
\le A\qw\sum_{k\ge1} k^2\knk
=A\qw S_2(F_n)
\\&
=nA\qw s_2(F_n),
  \end{split}
\end{equation*}
and thus the weight of any component $\cC$ in $H_n$ is
\begin{equation*}
  \begin{split}
\sum_{i\in\cC} x_i 
&\le 	
\sum_{i\in\cC} x_i \ett{x_i\le A}
+
\sum_{i} x_i \ett{x_i> A}
\le A|\cC|+nA\qw s_2(F_n)
\\&
\le A C_1(H_n)+nA\qw s_2(F_n).
  \end{split}
\end{equation*}
For any $\eps>0$, we may choose $A=A_n\=\eps\qw s_2(F_n)$ 
and find (since $A_n=O(1)$) \whp, using \eqref{c1small}, 
\begin{equation}\label{ms}
  C_1(\gynt)=\sup_{\cC}\sum_{i\in\cC} x_i 
\le
A_n C_1(H_n) + \eps n
\le 2\eps n.
\end{equation}
which proves (i) when $t\le1/\bs_2$.

\subsection{The supercritical case}

Suppose now that $t>\bs_2\qw$.

By \cite[Theorem 3.1]{SJ178},
the size $C_1(H_n)$ of the largest component $\cC_1$ of $H_n$ satisfies
\begin{equation*}
  \frac{C_1(H_n)}n\pto\rho(\kk)>0.
\end{equation*}
Furthermore $C_2(H_n)=\op(n)$, and it follows by the same argument as for
\eqref{ms} above that the weigth of any component $\cC\neq\cC_1$ of $H_n$ is
at most
\begin{equation*}
\max_{\cC\neq\cC_1}\sum_{i\in\cC} x_i 
\le
A_n C_2(H_n) + \eps n
\le 2\eps n
\end{equation*}
\whp, and thus $\op(n)$.
Since $\cC_1$ has weight $\ge|\cC_1|=\rho(\kk)n+\op(n)$, it follows that
\whp{} the largest component $\cC_1$ of $H_n$ also has the largest weight,
and thus corresponds to the largest component in $\gynt$, while 
$C_2(\gynt)=\op(n)$.

It remains to find the weight of $\cC_1$. We first note that by 
\cite[(2.13), Theorem 6.2 and (5.3)]{SJ178},
$\rho(\kk)=\int_\bbzpp\rhokk(x)\dd\mu(x)$, where $\rhokk(x)$ is the unique
positive solution to  
\begin{equation*}
  \rhokk=\Phi_\kk(\rhokk)
\=1-e^{-T_\kk\rhokk}.
\end{equation*}
Since 
\begin{equation*}
T_\kk\rhokk(x)\=\int_\bbzpp  \kk(x,y)\rhokk(y)\dd\mu(y)
=tx\int_\bbzpp  y\rhokk(y)\dd\mu(y),
\end{equation*}
we thus have 
\begin{equation*}
  \rhokk(x)=1-e^{-\rho t x}
\end{equation*}
with
\begin{equation}\label{ga}
  \rho=\int_\bbzpp x\rhokk(x)\dd\mu(x)
=\int_\bbzpp \rhokk(x)\dd\hmu(x)
=\int_\bbzpp \lrpar{1-e^{-\rho t x}}\dd\hmu(x).
\end{equation}

To find the weight $w(\cC_1)$ of $\cC_1(H_n)$, we note that if $f(x)\=x$, 
then $f:\bbzp\to\bbR$ satisfies, using \eqref{xdmu},
$\frac1n\sum_i f(x_i)=\frac1n\sum_i x_i =|F_n|/n=1=\int f\dd\mu$,
and thus \cite[Theorem 9.10]{SJ178} applies and yields 
\begin{equation}\label{wc1}
\frac1n  w(\cC_1)
=\frac1n\sum_{i\in\cC_1}x_i
\pto \int_\bbzpp x\rhokk(x)\dd\mu(x)
=\rho.
\end{equation}
Combining \eqref{ga} and \eqref{wc1}, we thus find that
\begin{equation}\label{emma}
  |C_1(\gynt)|=w(\cC_1(H_n))=\rho n+\op(n),
\end{equation}
where $\rho$ solves the equation \eqref{ga}, which also can be written
\begin{equation}\label{ga2}
  \rho=\E\lrpar{1-e^{-\rho t Z}}=1-\E e^{-\rho t Z}.
\end{equation}

Applying \refL{L1} to $Y\=t Z$, we see that when $t>1/\bs_2=1/\E Z$, there is a
unique $\rho>0$ satisfying \eqref{ga2}. 
%(Of course, $\rho=0$ is also a solution.)

Further, in \ref{t1gan}, 
we may apply \refL{L1} also to $Y\=t Z_n$; thus
there indeed is a
unique such $\rho_n$. Moreover, by \refL{L1}\ref{l1lim},
$\rho_n\to\rho$.
Hence, \eqref{emma} yields
\begin{equation*}
  |C_1(\gynt)|=\rho_n n+\op(n),
\end{equation*}
which proves \ref{t1gan} when $t>1/\bs_2$.

We have shown the conclusions in \ref{t1sub} and \ref{t1gan} when
$t\le1/\bs_2$ and $t>1/\bs_2$, respectively. However, the statements use
instead the slightly different conditions $t\le 1/s_2(F_n)$ and $t>1/s_2(F_n)$.
For \ref{t1sub}, this is no problem: if $t\le1/s_2(F_n)$ for infinitely many
$n$, then $t\le 1/\bs_2$ since we have assumed $s_2(F_n)\to\bs_2$.

To complete the proof of \ref{t1gan}, however, we have to consider also the
case $1/\bs_2\ge t> 1/s_2(F_n)$.
If this holds (for a subsequence), then
$\E(t Z_n)=ts_2(F_n)\le s_2(F_n)/\bs_2\to1$, and thus $\rho_n\to0$ by
\refL{L1}\ref{l1lim}. 
Since $t\le1/\bs_2$, \eqref{ms} applies and shows that
\begin{equation}\label{jesper}
  |C_1(\gynt)|=\op(n)=\rho_n n+\op(n),
\end{equation}
so \ref{t1gan} holds in this case too.
This completes the proof of \ref{t1sub} and \ref{t1gan}.

\ref{t1+-p} now follows easily from \refL{L1}.
We have, by \eqref{ezk}, $\E(tZ_n)=t s_2=1+\gd_ns_2$, $\E(tZ_n)^2=t^2s_3$ and
$\E(tZ_n)^3=t^3s_4$. Hence, 
\begin{equation}\label{erika}
  \frac{\E(tZ_n)-1}{\E(tZ_n)^2}
=\frac{\gd_n s_2}{t^2s_3} 
=\frac{\gd_n s_2^3}{(1+\gd_ns_2)^2s_3} 
> \gd_n\frac{ s_2^3}{s_3} (1-2\gd_ns_2), 
\end{equation}
so the lower bound follows by \ref{t1gan} and \refL{L1}\ref{l1-}.

For the upper bound we have by \eqref{erika}
\begin{equation*}
  \frac{\E(tZ_n)-1}{\E(tZ_n)^2}
%=\frac{\gd_n s_2^3}{(1+\gd_ns_2)^2s_3} 
< \gd_n\frac{ s_2^3}{s_3} , 
\end{equation*}
and similarly
\begin{equation*}
  \frac{(\E(tZ_n)-1)\E(tZ_n)^3}{(\E(tZ_n)^2)^2}
=\frac{\gd_n s_2 t^3s_4}{t^4s_3^2} 
=\frac{\gd_n s_2^2s_4}{(1+\gd_ns_2)s_3^2}
<\frac{\gd_n s_2^2s_4}{s_3^2},
\end{equation*}
and the upper bound follows by \refL{L1}\ref{l1+}.

For \ref{t1+-whp}, we note that 
if $\liminf_n\gd_n>0$, we can by ignoring some small $n$ assume that
$\inf_n\gd_n>0$, and then the difference between the \lhs{} and \rhs{} in
\eqref{erika} is bounded below (since $1\le s_2\le s_3=O(1)$); hence
we can add some small $\eta>0$ to the right hand side
of \eqref{erika} such that the inequality still holds for large $n$.
Consequently,
\begin{equation*}
C_1(\gynt)/n
  \ge \gd_n\frac{ s_2^3}{s_3} (1-2\gd_ns_2)-\eta+\op(1),
\end{equation*}
which implies that \whp
\begin{equation*}
C_1(\gynt)/n
  \ge \gd_n\frac{ s_2^3}{s_3} (1-2\gd_ns_2)
\end{equation*}
The upper bound follows in the same way.
\end{proof}

\section{Proof of Theorems \ref{TSa}--\ref{TBF}}\label{SpfBF}

\begin{proof}[Proof of \refT{TSa}]
  Define the functions
  \begin{align*}
f(t)&\=1/\bs_2(t),\\
g(t)&\=\bs_3(t)/\bs_2^3(t)=f^3(t)\bs_3(t),\\
h(t)&\=	\bs_4(t)/\bs_2^4(t)=f^4(t)\bs_4(t).
  \end{align*}
The differential equations \eqref{s2}--\eqref{s4} then translate into, after
simple calculations including some cancellations,
\begin{align}
  f'(t)&
%=-\frac{\bs_2'(t)}{\bs_2^2(t)}
=-\bx^2(t)f^2(t)-\bigpar{1-\bx^2(t)}, \label{f'}\\
  g'(t)&
%=3\bx^2(t)\bs_2\qwww(t)+3\bigpar{1-\bx^2(t)}\bs_2(t)\qww \bs_3(t)
=3\bx^2(t)f^3(t)-3\bx^2(t)f(t)g(t),     
\label{g'}\\ 
 h'(t)&=7\bx^2(t)f^4(t)+3\bigpar{1-\bx^2(t)}g^2(t)f\qww(t)-4\bx^2(t)f(t)h(t). 
\label{h'}
\end{align}

Consider first \eqref{f'}. 
The right hand side is locally Lipschitz in $t$ and $f$, and thus there
exists a unique solution with $f(0)=1$ in some maximal interval $[0,t_f)$ with
$t_f\le\infty$; if $t_f<\infty$ (which actually is the case, although we do
  not need this), $|f(t)|\to\infty$ as $t\upto t_f$.
Since $0<\bx(t)<1$ for all $t>0$, and further
$\bx(t)$ is decreasing, $f'(t)\le -(1-\bx^2(t))<-c_0$, for some $c_0>0$ and
all $t>0.1$, say. Hence, $f(t)$ decreases and will hit 0 at some finite time
$\tc<t_f$. This means that $\bs_2(t)=1/f(t)\to\infty$ as $t\upto\tc$, so
\eqref{s2} has a (unique) solution in $\otc$ but not further.

We have $f(\tc)=0$ and thus, by \eqref{f'},
$f'(\tc)=-(1-\bx^2(\tc))<0$. Consequently, defining $\rho\=(1-\bx^2(\tc))\qw>0$,
\begin{equation*}
 f(t)=\rho\qw(\tc-t)\bigpar{(1+O(\tc-t)},
\qquad t\le \tc, 
\end{equation*}
and thus
\begin{equation*}
 \bs_2(t)=\frac{\rho}{\tc-t}\bigpar{(1+O(\tc-t)},
\qquad t< \tc, 
\end{equation*}
as asserted.

Next, treating $\bx(t)$ and $f$ as  known functions, 
\eqref{g'} is a linear differential
equation in $g$. An integrating factor is
\begin{equation}\label{G}
  G(t)\=3\int_0^t  \bx^2(u)f(u)\dd u,
\end{equation}
and then the unique solution in $\otf$ is given by
\begin{equation}\label{g}
  g(t)=e^{-G(t)}+
3e^{-G(t)}\int_0^t e^{G(u)} \bx^2(u)f^3(u)\dd u.
\end{equation}
Hence \eqref{s3} has the unique solution $g(t)\bs_2^3(t)$, $t\in\otc$,
with $g(t)$ given by \eqref{g}. Note that $g(t)>0$ for $t\le\tc$.

Let $\gb\=g(\tc)>0$.
By \eqref{g'}, $g'(\tc)=0$, and thus, for $t<\tc$,
$g(t)=\gb+O(\tc-t)^2$, and 
\begin{equation*}
  \bs_3(t)=\gb \bs_2^3(t)\bigpar{1+O(\tc-t)^2}
=\frac{\gb\ga^3}{(\tc-t)^3}\bigpar{(1+O(\tc-t)},
\end{equation*}

Finally we consider \eqref{h'}. Here the \rhs{} is singular at $\tc$ because
of the factor $f\qww(t)$ in the second term, so we modify $h$ and consider
\begin{equation*}
h_1(t)\=h(t)-3g^2(t)\bs_2(t)  
=h(t)-3g^2(t)f\qw(t),
\end{equation*}
which satisfies the differential equation
\begin{equation*}
  \begin{split}
h_1'(t)
&=7\bx^2(t)f^4(t)	-18\bx^2(t)g(t)f^2(t)+15\bx^2(t)g^2(t)		
-4\bx^2(t)f(t)h(t)	
\\&
=7\bx^2(t)f^4(t)	-18\bx^2(t)g(t)f^2(t)+3\bx^2(t)g^2(t)		
-4\bx^2(t)f(t)h_1(t).	
  \end{split}
\end{equation*}
Again, this is a linear differential equation, with a unique solution in 
$\otf$. We leave the explicit form to the reader, since we need only that
$h_1(t)=O(1)$ for $t\le\tc$, which yields that for $t\in\otc$,
\begin{equation*}
  \begin{split}
  \bs_4(t)&=h(t)\bs_2^4(t)=3g^2(t)\bs_2^5(t)+h_1(t)\bs_2^4(t)	
\\&=3\gb^2 \bs_2^5(t)+O\bigpar{\bs_2^4(t)}.
\qedhere
  \end{split}
\end{equation*}
\end{proof}

\begin{comment}
 given by, noting that $h_1(0)=h(0)-3=-2$,
\begin{multline*}
  h_1(t)=-2e^{-H(t)}+e^{-H(t)}
\cdot\\
\intot e^{H(u)}
\Bigpar{7\bx^2(u)f^4(u)	-18\bx^2(u)g(u)f^2(u)+3\bx^2(u)g^2(u)}\dd u
\end{multline*}
where
\begin{equation*}
  H(t)\=4\intot \bx^2(t)f(t) \dd t.	
\end{equation*}
\end{comment}

\begin{proof}[Proof of \refT{TSb}]
  For $k=2$, this is, as  said above, proved in 
\cite[Theorems 1.1 and 4.3]{SW}.
We prove the extension by the same method (with somewhat different
notation).

Let, for a vertex $v\in G$,  
$\cC(v)$ be the component of $G$ containing the vertex $v$, and
$C(v)\=|\cC(v)|$. 

For a given graph $G$, let $G^+$ be the random graph obtained by adding one
random edge by the \BFq{} rule; we assume that the edge was chosen from the
pair $e_1=\set{v_1,w_1}$ and $e_2=\set{v_2,w_2}$. 
If the added edge is $\set{v,w}$ (which thus is either $\set{v_1,w_1}$ or
$\set{v_2,w_2}$), and further $\cC(v)\neq\cC(w)$, then, by \eqref{Sk}, 
\begin{equation}\label{sk+}
  S_k(G^+)-S_k(G)=\bigpar{C(v)+C(w)}^k-C(v)^k-C(w)^k,
\end{equation}
while $S_k(G^+)-S_k(G)=0$ if $\cC(v)=\cC(w)$.
We define
\begin{equation}\label{gdx}
  \gdx_k=\gdx_k(G;v,w)\=\bigpar{C(v)+C(w)}^k-C(v)^k-C(w)^k.
\end{equation}
Hence,
\begin{equation*}
  \begin{split}
  \E\bigpar{S_k(G^+)-S_k(G)-\gdx_k}
&
=-\E\bigpar{\gdx_k\ett{\cC(v)=\cC(w)}}
\\&
=-\E\bigpar{(2^k-2)C(v)^k\ett{\cC(v)=\cC(w)}}	
  \end{split}
\end{equation*}
and thus
\begin{equation*}%\label{ba2}
  \begin{split}
| \E&\bigpar{S_k(G^+)-S_k(G)-\gdx_k}|
\le2^k\E\bigpar{C(v)^k\ett{\cC(v)=\cC(w)}}	
\\&
\le2^k\E\bigpar{C(v_1)^k\ett{\cC(v_1)=\cC(w_1)}}	
+2^k\E\bigpar{C(v_2)^k\ett{\cC(v_2)=\cC(w_2)}}	
\\&
=\frac{2^{k+1}}n\E C(v_1)^{k+1}
%=\frac{2^{k+1}}{n^2}S_{k+2}(G)
%=O\lrparfrac{C_1(G)^{k+2}}{n}
\le 2^{k+1}\frac{C_1(G)^{k+1}}{n}
.
  \end{split}
\end{equation*}
In particular, if $C_1(G)=O(\log n)$, then
%$S_{k+1}(G)=O\bigpar{n\log^{k+1}n}$ and 
\begin{equation}\label{ba2b}
  \begin{split}
| \E&\bigpar{S_k(G^+)-S_k(G)-\gdx_k}|
=O\Bigparfrac{\log^{k+1} n}{n}
=o(1).
  \end{split}
\end{equation}

Expanding \eqref{gdx}, we have
\begin{align}
  \gdx_2&=2C(v)C(w),\label{ba2a2}\\
  \gdx_3&=3C(v)^2C(w)+3C(v)C(w)^2,\label{ba2a3}\\
  \gdx_4&=4C(v)^3C(w)+6C(v)^2C(w)^2+4C(v)C(w)^3.
\label{ba2a4}
\end{align}

The \BFq{} rule is to take $\set{v,w}=\set{v_1,w_1}$ if $C(v_1)=C(w_1)=1$.
The probability of this is $x_1(G)^2$, and in this case $\gdx_k=2^k-2$.

The opposite case $\set{v,w}=\set{v_2,w_2}$, which we denote by $\cE_2$, 
has probability $1-x_1(G)^2$. 
Conditioning on this case places us basically
in the well-studied Erd\H{o}s--R\'enyi regime.
That is, $v$ and $w$ are uniform and independent,
and thus for any $k$ and $\ell$,
\begin{equation*}
  \begin{split}
  \E\bigpar{C(v)^kC(w)^l\mid\cE_2}
&=\frac1{n^2}\sum_{v,w}C(v)^kC(w)^\ell
=\frac1{n^2}\sum_{i}C_i^{k+1}\sum_{j}C_j^{\ell+1}
\\&
=s_{k+1}(G)s_{\ell+1}(G).	
  \end{split}
\end{equation*}
Hence, \eqref{ba2a2}--\eqref{ba2a4} yield
\begin{align*}
  \E\gdx_2&=2x_1^2(G)+(1-x_1(G)^2)\cdot2s_2(G)^2 ,%\label{ba2aa2} 
\\
  \E\gdx_3&=6x_1^2(G)+(1-x_1(G)^2)\cdot6s_2(G)s_3(G),%\label{ba2aa3} 
\\
\E\gdx_4&=14x_1^2(G)+(1-x_1(G)^2)\cdot\bigpar{8s_2(G)s_4(G)+6s_3(G)^2}.
%\label{ba2aa4} 
\end{align*}
By \eqref{ba2b}, we thus have, for $k=2,3,4$ and provided $C_1(G)=O(\log n)$,
\begin{multline}\label{ba3}
%  \begin{split}
\E\bigpar{S_k(G^+)-S_k(G)}
=\E\gdx_k+O\bigpar{\log^{k+1} n/n}
\\
=2f_k\bigpar{x_1(G),s_2(G),s_3(G),s_4(G)}+O\bigpar{\log^{k+1} n/n},
%  \end{split}
\end{multline}
%(uniformly in all $G$ with $|G|=n$),
with
\begin{align*}
f_2(x_1,s_2,s_3,s_4)&\=x_1^2+(1-x_1^2)s_2^2 ,\\
f_3(x_1,s_2,s_3,s_4)&\=3x_1^2+3(1-x_1^2)s_2s_3,\\
f_4(x_1,s_2,s_3,s_4)&\=7x_1^2+(1-x_1^2)\bigpar{4s_2s_4+3s_3^2} .
\end{align*}
Similarly, as shown in \cite{SW},
\begin{equation}
\label{ba3b}
\E\bigpar{n_1(G^+)-n_1(G)}
=2f_1\bigpar{x_1(G),s_2(G),s_3(G),s_4(G)}+O\bigpar{1/n},
\end{equation}
where (the variables $s_2,s_3,s_4$ are redundant here)
\begin{align*}
  f_1(x_1,s_2,s_3,s_4)&\=-x_1^2-(1-x_1^2)x_1 .
\end{align*}

Consider the vector-valued random process 
$$
X_i\=\bigpar{x_1(\BF_i),s_2(\BF_i), s_3(\BF_i), s_4(\BF_i)},
$$
and let $\cF_i=\gs(X_0,\dots,X_i)$ be the $\gs$-field describing the history
up to time $i$. Further, let $\ff\=(f_1,f_2,f_3,f_4):\bbR^4\to\bbR^4$.
Using this notation, \eqref{ba3}--\eqref{ba3b} yield
\begin{equation}
\label{ba3c}
\E\bigpar{n(X_{i+1}-X_i)\mid\cF_i}-2\ff(X_i)
=O\bigpar{\log^5 n/n},
\end{equation}
uniformly in $i\le tn/2$, provided $C_1(\BF_i)=O(\log n)$.

By \cite[Theorem 1.1]{SW}, there exists a constant $c'$ (depending on $t$)
such that \whp{} $C_1(\BF_i)\le c'\log n$ for all $i\le tn/2$.
As in \cite{SW}, we avoid the problem when $C_1(\BF_i)>c'\log n$ by defining
$\XX_0=X_0=(1,1,1,1)$, $\XX_{i+1}=X_{i+1}$ when $C_1(\BF_i)\le c'\log n$
and $\XX_{i+1}=\XX_i+\frac2n \ff(\XX_i)$ otherwise.
Then \whp{} $\XX_i=X_i$ for all $i\le tn/2$, so we can just as well consider
$\XX_i$. We have, by \eqref{ba3c} but now without side condition,
for all $i\le tn/2$,
\begin{equation*}
  \E\bigpar{n(\XX_{i+1}-\XX_i)\mid\cF_i}
=2 \ff(\XX_i)+O(\log^5n/n)
\end{equation*}
and also, for some $c''$, from \eqref{sk+} and $|n_1(G^+)-n_1(G)|\le2$,
\begin{equation*}
  |\XX_{i+1}-\XX_i|\le c''\log^4n/n.
\end{equation*}
The differential equation method in the form of \citetq{Theorem 4.1}{SW},
which is taken from \citetq{Theorem 5.1}{W99}, now applies (with
$Y(i)=n\XX_i$)
and the result follows; note that the differential equations
\eqref{x1}--\eqref{s4} can be written
$\gf'(t)=\ff(\gf(t))$ with $\gf=(\bx,\bs_2,\bs_3,\bs_4)$,
where further $\gf(0)=(1,1,1,1)=X_0=\XX_0$.
\end{proof}

\begin{proof}[Proof of \refT{TBF}]
We may assume that $\gd$ is small, since the result is trivial for
$\gd\ge\gd_0>0$ if we choose $\C$ large enough. In particular, we assume
$\gd<1$. 

Let $\eps\=\gd^{2/3}>\gd$.
We stop the process at $\tc-\eps$, and let $F:=\BF(\tc-\eps)$. We then let 
the process evolve to $\tc+\gd$ by adding $(\eps+\gd)n/2$ further edges
according to the Bohman--Frieze rule. Actually, for convenience, we add
instead a random number of edges with a Poisson distribution
$\Po\bigpar{(\eps+\gd)n/2}$; this will not affect our asymptotic results (by
the same standard argument as for comparing the different models in
\refS{SER}).
We denote the resulting graph by $\tbfd$.

By Theorems \ref{TSb} and \ref{TSa}, for $k=2,3,4$,
and with $\cg_k$ as in \refT{TSa},
\begin{equation*}
  s_k(F)=\bs_k(\tc-\eps)+\op(1)=\frac{\cg_k}{\eps^{2k-3}}\bigpar{1+O(\eps)}
+\op(1).
\end{equation*}
Since $|\op(1)|\le\eps$ \whp, we thus have \whp{}  
\begin{equation}\label{ca1}
  s_k(F)=\frac{\cg_k}{\eps^{2k-3}}\bigpar{1+O(\eps)}.
\end{equation}
(This means that there exists a constant $c$, not depending on $\eps$ or
$n$, such that \eqref{ca1} holds with the error term
$O(\eps)\in[-c\eps,c\eps]$ w.h.p.)
Similarly, $x_1(F)\pto\bx(\tc-\eps)=\bx(\tc)+O(\eps)$,
so \whp{} $x_1(F)=\bx(\tc)+O(\eps)$.

We fix $F$ (\ie, we condition on $F$) and assume that \eqref{ca1} holds
together with $x_1(F)=\bx(\tc)+O(\eps)$ (for some fixed  implicit
constant $c$ in the $O(\eps)$; we have just shown that this holds \whp\
provided $c$ is chosen large enough).

We cannot directly apply \refT{T1} since the graph evolves by the \BFq{}
evolution and not by the \ER{} evolution. Nevertheless, we can approximate
and find upper and lower bounds of the graphs
where we can apply \refT{T1}; the idea is
that we consider the \ER{} edges separately as an \ER{} evolution.

For a lower bound, let $V_1$ be the set of isolated vertices in $F$ and
consider only the pairs of edges 
$e_1=\{v_1,w_1\}$, $e_2=\{v_2,w_2\}$ where $v_1\notin V_1$ or $w_1\notin
V_1$.
Since the graphs $\BF_\ell$ in the continued process contain $F$, the vertices 
$v_1$ and $w_1$ are not both isolated in the current $\BF_\ell$,
and thus $e_2=(v_2,w_2)$ is added, and these are independent \ER{} edges,
\ie, uniformly chosen.
The number of such \ER{} edges is $\Po\bigpar{(1-x_1(F)^2)\eps+\gd)n/2}$,
since each time we add an edge, the probability of it being of this type 
is $1-(|V_1|/n)^2=1-x_1(F)^2$. (Note that we ignore some \ER{} edges in
order to avoid unpleasant dependencies.)

Call the resulting graph $\BFG^-\subseteq\BF(\tc+\gd)$. Then
\refT{T1}\ref{t1+-whp} applies to $\BFG^-$, with
\begin{equation*}
  t=\bigpar{1-x_1(F)^2}(\eps+\gd)
=\bigpar{1-\bx(\tc)^2+O(\eps)}(\eps+\gd)
\end{equation*}
and, recalling \eqref{ca1} and $\ga=\bigpar{1-\bx^2(\tc)}\qw$,
\begin{equation}\label{gd-}
  \gd_n=t-1/s_2(F)
=\bigpar{1-\bx(\tc)^2}(\eps+\gd)-\ga\qw\eps+O(\eps^2)
=\bigpar{1-\bx(\tc)^2}\gd+O(\eps^2),
\end{equation}
which yields \whp, using again \eqref{ca1},
\begin{equation}\label{c1low}
  \begin{split}
\frac{C_1(\tbfd)}n
&\ge
\frac{C_1(\BFG^-)}n
\ge 2\gdn\frac{s_2(F)^3}{s_3(F)}\bigpar{1-2\gdn s_2(F)} %-\op(1)
\\&
=2\lrpar{\bigpar{1-\bx(\tc)^2}\gd+O(\eps^2)}
 \frac{\ga^3\eps^{-3}}{\gb\ga^3\eps^{-3}}
 \lrpar{1+O(\eps)+O\parfrac{\gd+\eps^2}{\eps}}
\\&
=\frac2{\ga\gb}\gd\lrpar{1+O(\eps)+O(\gd/\eps)+O(\eps^2/\gd)}
\\&
=\frac2{\ga\gb}\gd\lrpar{1+O(\gd\qqq)}
=\gam\gd+O(\gd\qqqd),
  \end{split}
\raisetag{\baselineskip}
\end{equation}
with our choice $\eps=\gd^{2/3}$ (which is optimal in this estimate).
 
For an upper bound, note that \whp{} at most $(\eps+\gd)n\le2\eps n$ edges
are added to $F$, so at most $4\eps n$ vertices are hit,
and thus during the process from $F$ to $\tbfd$,
$$
x_1\ge x_1(F)-4\eps=\bx(\tc)-O(\eps).
$$
Hence we add \whp{} at most 
$$
\bigpar{1-(\bx(\tc)-O(\eps))^2}(\eps+\gd)n/2
=
\bigpar{1-\bx(\tc)^2+O(\eps)}(\eps+\gd)n/2
$$
\ER{} edges. 
We also add a number of non-\ER{} edges, all joining two isolated vertices
(or being loops). They may depend on the \ER{} edges already chosen, but we
avoid this dependency by being generous and adding the edge $e_1=(v_1,w_1)$
in each round whenever both $v_1$ and $w_1$ are isolated in $F$ and neither
is an endpoint of an already added non-\ER{} edge. (We add $e_2$ by the same
\BFq{} rule as before, so we may now sometimes add both $e_1$ and $e_2$.)

%For convenience we add some extra edges, and 
Let $c_1$ be a large constant and let $\BFG^+$ be the graph obtained from $F$
by adding  
$2\eps n$ (to be on the safe side)
non-\ER{} edges in this way, together with
$\bigpar{1-\bx(\tc)^2+c_1\eps}(\eps+\gd)n/2$ \ER{} edges, independent of
each other and of the non-\ER{} edges.
We conclude that, if $c_1$ is chosen large enough,
we may couple $\BFG^+$ with the \BFp{} such that 
\whp{} $\tbfd\subseteq \BFG^+$.

Since the two types of edges are added independently, we may further add all
non-\ER{} edges first. Let $F_1$ be $F$ together with all non-\ER{} edges.
There are  $2\eps n$ such edges, and each joins two
isolated vertices and changes $S_k$ by $2^k-2$ (or by 0 if the edge is a
loop).
Hence,  for every $k\le4$, by \eqref{ca1},
\begin{equation}\label{skf1}
  s_k(F_1)=s_k(F)+O(\eps)
=\frac{\cg_k}{\eps^{2k-3}}\bigpar{1+O(\eps)}.
\end{equation}

Since $\BFG^+$ is obtained by adding the \ER{} edges to $F_1$, \refT{T1}
applies with 
\begin{equation*}
  t=\bigpar{1-\bx(\tc)^2+c_1\eps}(\eps+\gd)
\end{equation*}
and
\begin{equation}\label{gd+}
\gd_n=t-1/s_2(F_1)
=\bigpar{1-\bx(\tc)^2}(\eps+\gd)-\ga\qw\eps+O(\eps^2)
=\bigpar{1-\bx(\tc)^2}\gd+O(\eps^2),
\end{equation}
the same estimate as was obtained in \eqref{gd-}.
We use the upper bound in \refT{T1}\ref{t1+-whp}. 
By \eqref{skf1} and \eqref{gd+},
\begin{equation*}
\gdn  \frac{s_2(F_1)^2s_4(F_1)}{s_3(F_1)^2}
=  \gdn \frac{\cg_2^2\cg_4/\eps^7}{\cg_3^2/\eps^6}\bigpar{1+O(\eps)}
=O\Bigparfrac{\gdn}{\eps}
=O\Bigparfrac{\gd}{\eps}
=O\bigpar{\gd\qqq}.
\end{equation*}
Hence \refT{T1}\ref{t1+-whp} applies (for small $\gd$) and
yields, \whp, 
\begin{equation*}
  \begin{split}
n\qw C_1(\tbfd) 
&\le n\qw C_1(\BFG^+)	
\le 2\gdn \frac{s_2(F_1)^3}{s_3(F_1)}\bigpar{1+O(\gd\qqq)}
\\&
=2(1-\bx(\tc)^2)\gd \frac{\ga^3\eps^{-3}}{\gb\ga^3\eps^{-3}}
\bigpar{1+O(\eps^2/\gd+\eps+\gd\qqq)}
\\&
=\frac{2(1-\bx(\tc)^2)}\gb\gd \bigpar{1+O(\gd\qqq)}
=\gam\gd+O(\gd\qqqd).
  \end{split}
\end{equation*}
This and the corresponding lower bound \eqref{c1low} yield the result.
\end{proof}

\newcommand\AAP{\emph{Adv. Appl. Probab.} }
\newcommand\JAP{\emph{J. Appl. Probab.} }
\newcommand\JAMS{\emph{J. \AMS} }
\newcommand\MAMS{\emph{Memoirs \AMS} }
\newcommand\PAMS{\emph{Proc. \AMS} }
\newcommand\TAMS{\emph{Trans. \AMS} }
\newcommand\AnnMS{\emph{Ann. Math. Statist.} }
\newcommand\AnnPr{\emph{Ann. Probab.} }
\newcommand\CPC{\emph{Combin. Probab. Comput.} }
\newcommand\JMAA{\emph{J. Math. Anal. Appl.} }
\newcommand\RSA{\emph{Random Struct. Alg.} }
\newcommand\ZW{\emph{Z. Wahrsch. Verw. Gebiete} }
\newcommand\DMTCS{\jour{Discr. Math. Theor. Comput. Sci.} }

\newcommand\AMS{Amer. Math. Soc.}
\newcommand\Springer{Springer-Verlag}
\newcommand\Wiley{Wiley}

\newcommand\vol{\textbf}
\newcommand\jour{\emph}
\newcommand\book{\emph}
\newcommand\inbook{\emph}
\def\no#1#2,{\unskip#2, no. #1,} %(typeset after year) 
\newcommand\toappear{\unskip, to appear}

\newcommand\webcite[1]{%\hfil  %???
   %\penalty0 %???
\texttt{\def~{{\tiny$\sim$}}#1}\hfill\hfill}
\newcommand\webcitesvante{\webcite{http://www.math.uu.se/~svante/papers/}}
\newcommand\arxiv[1]{\webcite{arXiv:#1.}}

\def\nobibitem#1\par{}


\begin{thebibliography}{99}

\bibitem[Achlioptas, D'Souza and Spencer (2009)]{SSA}
D.~Achlioptas, R.~D'Souza and J.~Spencer,
Explosive percolation in random networks. 
\emph{Science} \vol{323} (2009), no. 5920, 1453--1455.


\bibitem[Bohman and Frieze(2001)]{BF}
T. Bohman and A. Frieze, 
Avoiding a giant component.  
\emph{Random Structures Algorithms}  \vol{19}  (2001),  no. 1, 75--85. 

\bibitem[Bohman, Frieze and Wormald(2004)]{BFW}
T. Bohman, A. Frieze and N.C. Wormald,
Avoidance of a giant component in half the edge set of a random graph.
\emph{Random Structures Algorithms}  \vol{25} (2004), no. 4, 432--449. 


\bibitem[Bollob\'as(2001)]{Bollobas}
B. Bollob\'as, 
\book{Random Graphs}, 2nd ed., Cambridge Univ. Press,
Cambridge, 2001.

\bibitem[Bollob\'as, Janson and Riordan(2007)]{SJ178}
B.~Bollob\'as, S. Janson and O.~Riordan,
The phase transition in inhomogeneous random graphs.
\RSA \vol{31} (2007), 3--122.

\bibitem[Erd\H os and R\'enyi(1959)]{ER59} 
P.~Erd\H os \& A.~R\'enyi,
On random graphs. I.  \emph{Publ. Math. Debrecen}  \vol6  (1959), 290--297. 


\bibitem[Erd\H os and R\'enyi(1960)]{ER60} 
P.~Erd\H os \& A.~R\'enyi,
 On the evolution of random graphs,
 \jour{Magyar Tud. Akad. Mat. Kutat\'o Int. K\"ozl.}
 \textbf{5} (1960), 17--61.

\nobibitem[Erd\H os and R\'enyi(1961)]{ER61}
P.~Erd\H os \& A.~R\'enyi,
 On the evolution of random graphs,
 \jour{Bull. Inst. Internat. Statist.}
 \textbf{38} (1961), 343--347.

\bibitem[Gut (2005)]{Gut}
A. Gut,
\book{Probability: A Graduate Course}.
Springer, New York, 2005.

\bibitem[Janson(2009+)]{SJN6}
S. Janson,
Probability asymptotics: notes on notation.
Institut Mittag-Leffler preprint 31, 2009 spring.

\bibitem[Janson(2009+)]{SJ241}
S. Janson,
Susceptibility of random graphs with given vertex degrees,
\arxiv{0911.2636}

\bibitem[Janson and Luczak(2008)]{SJ218}
S. Janson \& M. Luczak,
Susceptibility in subcritical random graphs.
\jour{J. Math. Phys.} \vol{49}:12 (2008), 125207.  

\bibitem[Janson, {\L}uczak and Ruci\'nski(2000)]{JLR}
S. Janson, T. \L uczak \& A. Ruci\'nski,
\book{Random Graphs}.
\Wiley, New York, 2000.

\bibitem[Janson and Riordan(2009+)]{SJ232}
S. Janson \& O.~Riordan,
Susceptibility in inhomogeneous random graphs.
\arxiv{0905.0437}


\nobibitem{Kallenberg}
O. Kallenberg,
\book{Foundations of Modern Probability.}
2nd ed., Springer, New York, 2002. 

\nobibitem{KnuthI} 
D.E. Knuth, 
\emph{The Art of Computer Programming. Vol. 1:
 Fundamental Algorithms}. 
3nd ed., Addison-Wesley,
Reading, Mass., 1997.

\bibitem[Spencer and Wormald(2007)]{SW}
J. Spencer and N. Wormald,
Birth control for giants.  
\emph{Combinatorica}
\vol{27}  (2007),  no. 5, 587--628. 

\bibitem[Wormald(1999)]{W99}
N. Wormald,
The differential equation method for random graph processes and greedy
algorithms.
In \emph{Lectures on Approximation and Randomized Algorithms},
M. Karo\'nski and H.-J. Pr\"omel, eds, pp. 73--155, PWN, Warsaw, 1999. 


\end{thebibliography}
\end{document}